\documentclass[12pt]{amsart}

\usepackage{amssymb, latexsym}

\newtheorem{Theorem}{Theorem}

\newtheorem{Lemma}[Theorem]{Lemma}
\newtheorem{Proposition}[Theorem]{Proposition}

\usepackage{hyperref}
\usepackage{graphicx}

\begin{document}

\title[Symmetry of zeros of Lerch zeta-function ] {Symmetry of zeros of Lerch zeta-function for equal parameters}

\author{ Ram\={u}nas Garunk\v{s}tis}
\address{Ram\={u}nas Garunk\v{s}tis \\
Department of Mathematics and Informatics, Vilnius University \\
Naugarduko 24, 03225 Vilnius, Lithuania}
\thanks{The first author is  supported
by grant No. MIP-049/2014 from the Research 
Council of Lithuania.}
\email{ramunas.garunkstis@mif.vu.lt}
\urladdr{www.mif.vu.lt/~garunkstis}

\author{ Rokas Tamo\v si\=unas}
\address{ Rokas Tamo\v si\=unas \\
Department of Mathematics and Informatics, Vilnius University \\
Naugarduko 24, 03225 Vilnius, Lithuania}
\email{rokas.tamosiunas@mif.vu.lt}



\maketitle

\begin{abstract} 
For most values of parameters $\lambda$ and $\alpha$, the zeros of the Lerch zeta-function $L(\lambda, \alpha, s)$ are distributed very chaotically.  In this paper we consider the special case of equal parameters $L(\lambda, \lambda, s)$  and show by calculations that the nontrivial zeros  either lie extremely close to the critical line $\sigma = 1/2$
or are distributed almost symmetrically with respect to the critical line. We also investigate this phenomenon theoretically.
\end{abstract}

\section{Introduction}

Let $s=\sigma+it$. Denote by $\{\lambda\}$ the fractional part of a real number $\lambda$. In this paper, $\varepsilon$ is any positive real number and $T$ always tends to plus infinity. In all theorems and lemmas, the numbers $\lambda$ and $\alpha$ are fixed constants.

For $0<\lambda, \alpha\leq 1$, the Lerch zeta-function is given by
\begin{align*}
L(\lambda,\alpha,s)=\sum_{m=0}^\infty\frac{e^{2\pi i\lambda m}}{(m+\alpha)^s} \qquad(\sigma>1).
\end{align*}
This function has analytic continuation to the whole complex plane except for a possible simple pole at $s=1$ (Lerch \cite{Lerch1887}, Laurin\v cikas and Garunk\v stis  \cite{gl}). 

 Let $l$ be a straight line in the complex plane ${\Bbb C}$,
and denote by $\varrho(s,l)$ the distance of $s$ from $l$.
Define, for $\delta>0$,
\begin{eqnarray*}
L_\delta(l)=\big\{s\in{\Bbb C}:\;\varrho(s,l)<\delta\big\}.
\end{eqnarray*}
In Garunk\v stis and Laurin\v cikas \cite{gar}, Garunk\v stis and Steuding \cite{Garunkstis2002}, for $0<\lambda<1$ and $\lambda\ne1/2$,   it is proved that 
$L(\lambda,\alpha,s)\ne0$ if $\sigma<-1$ and
\begin{eqnarray*}
s\not\in L_{\log4\over\pi}\bigg(\sigma =\frac{\pi
t}{\log{\frac{1-\lambda}{\lambda}}}+1\bigg)\,.
\end{eqnarray*}
  For $\lambda=1/2,1$, from Spira \cite{Spira} and \cite{gar} we see that  $L(\lambda,\alpha,s)\ne0$ if $\sigma<-1$ and $|t|\ge1$.  Moreover, in  \cite{gar} it is showed that $L(\lambda,\alpha, s)\ne0$ if $\sigma\ge1+\alpha$.  We say that a zero of $L(\lambda,\alpha, s)$ is {\it nontrivial} if it lies in the strip $-1\le\sigma<1+\alpha$ and we denote a nontrivial zero by $\rho=\beta+i\gamma$.

Let $\zeta(s)$ and $L(s,\chi)$ denote the Riemann zeta-function and the Dirichlet $L$-function accordingly. We have that 
\begin{align}\label{RH}
L(1,1,s)=\zeta(s)\quad \text{and}\quad L(1/2,1/2,s)=2^s L(s,\chi),
\end{align} 
where $\chi$ is a Dirichlet character$\mod 4$
with $\chi(3)=-1$.  For these two cases, certain versions of the Riemann hypothesis (RH) can be formulated. Similar cases are $L(1,1/2,s)=(2^s-1)\zeta(s)$ and  $L(1/2,1,s)=(1-2^{1-s})\zeta(s)$.  
For all the other cases,  it is expected that the real parts of zeros of the Lerch zeta-function form a dense subset of the interval $(1/2,1)$. This is proved for any $\lambda$ and transcendental $\alpha$ (\cite[Theorem 4.7 in Chapter 8]{gl}).

In this paper, we investigate the zero distribution of the Lerch zeta-function $L(\lambda, \alpha, s)$ when the parameters are equal, i.e. $\lambda=\alpha$. The motivation for this are calculations which show that the first nontrivial zeros of $L(\lambda, \lambda, s)$  are often located  almost on the critical line $\sigma=1/2$. Next are the first 4 zeros (rounded to two decimal numbers) of several  Lerch zeta-functions.

$L(1/3,1/3, s) :$ $0.50+3.99 i$, $0.50+7.28 i$, $0.50+9.54 i$, $0.50+12.18 i$.

$L(1/3,2/3, s) :$ $0.86+5.68 i$, $0.53+9.59 i$, $0.86+12.66 i$,  $0.49+15.11 i$.

$L(3/4,3/4, s) :$ $0.50+9.69 i$, $0.50+15.26 i$, $0.50+18.65 i$, $0.50+23.05 i$.

$L(1/4,3/4, s) :$ $1.03+5.24 i$, $0.64+8.81 i$, $0.76+11.96 i$,  $0.88+14.19 i$.

For  a rational number $\lambda\ne1/2,1$ it is expected that the function $L(\lambda,\lambda, s)$ has many zeros off the critical line. Our calculations then show that the zeros are  almost symmetrically distributed with respect to the critical line. For example, for $L(3/4,3/4, s)$ we have the following zeros: $-0.10+120.60i$ and $1+0.10+120.60i$; $0.37+202.77i$ and $1-0.37+202.77i$. Usually such symmetry of zeros can be explained by the shape of the functional equation. A typical example is the Heillbronn Davenport zeta-function. Possibly such symmetry  forces zeros to stay  on the critical line more often. More on this see, for example, Bombieri and Hejhal \cite{bh}, Balanzario and S\'anchez-Ortiz \cite{Balanzario2007}, Garunk\v stis and \v Sim\.enas \cite{gs}, Vaughan \cite{vaug}. For the Lerch zeta-function the following relation, usually called the functional equation, is true.
\begin{align}\label{Lerchfunc}
L(\lambda,\alpha,1-s)=&(2\pi)^{-s}\Gamma(s)\biggr(
e^{\pi i\frac{s}{2}-2\pi i\alpha\lambda}L(1-\alpha,\lambda,s)
\\&
+e^{-2\pi i\frac{s}{4}+2\pi i \alpha(1-\{\lambda\})}L(\alpha,1-\{\lambda\},s)\biggr).\nonumber
\end{align}
Various proofs of this functional equation can be found in Lerch \cite{Lerch1887}, Apostol \cite{Apostol1951}, Oberhettinger \cite{Oberhettinger1956}, Mikol\'as \cite{Mikolas1971}, Berndt \cite{Berndt1972}, see also Lagarias and Li \cite{ll1}, \cite{ll2}. The Lerch zeta function has a second moment (Garunk\v{s}tis, Laurin\v{c}ikas, and Steuding \cite{gls2003}) and it is a universal function (Laurin\v{c}ikas \cite{lau}, Lee,  Nakamura,  Pa\'nkowski \cite{lnp}).

  For $\lambda=\alpha$, we can rewrite (\ref{Lerchfunc}) as
\begin{align}\label{almostsymetry}
\overline{L(\lambda,\lambda,1-\overline{s})}=&(2\pi)^{-s}\Gamma(s)e^{-\pi i\frac{s}{2}+2\pi i \lambda^2}L(\lambda,\lambda,s)\nonumber
\\&
+(2\pi)^{-s}\Gamma(s)
e^{\pi i\frac{s}{2}-2\pi i(1-\lambda)\lambda}L(1-\lambda,1-\{\lambda\},s)
\\
=&G(s)L(\lambda,\lambda,s)+P(s).\nonumber
\end{align}
By the bound for the Lerch zeta-function and by Stirling's formula we see that, for any vertical strip, $|P(s)|< t^Ae^{-\pi t}$ and $|G(s)|\ge t^B$ (see Lemma \ref{lerchgrow} and its proof below). Thus the shape of the formula (\ref{almostsymetry})  suggests that the nontrivial zeros of $L(\lambda,\lambda,s)$ should be distributed almost symmetrically with the respect of the critical line. However calculations in the next section show that this symmetry is not strict.

Denote by $N(\lambda,\alpha,T)$  the number of 
nontrivial zeros of the function $L(\lambda,\alpha,s)$ in the region 
$0<t<T$. 
For $0<\lambda,\alpha\leq 1$,  we have (\cite{gar})
\begin{align}\label{zeronumber}
N(\lambda,\alpha,T)=\frac{T}{2\pi}\log 
\frac{T}{2\pi e\alpha\lambda}+O(\log T).
\end{align}

 The next theorem shows that in the upper half-plane nontrivial zeros of the Lerch zeta-function with equal parameters  on average are symmetrically distributed with a small error term.
\begin{Theorem}\label{averagesymm}
For $0<\lambda,\alpha\leq 1$, 
\begin{eqnarray*}
\sum_{0<\gamma\leq T}\left(\beta-{1\over 2}\right)
={T\over 4\pi}\log{\alpha\over \lambda}+O(\log T).
\end{eqnarray*}
\end{Theorem}

Now we consider the symmetry of  the individual zeros. Let $\rho$ be a zero of $L(\lambda, \lambda, s)$. In view of (\ref{almostsymetry}) and Rouch\'e's theorem we see that $L(\lambda, \lambda, s)$ has an almost symmetrical zero in some small disc $|s-(1-\overline{\rho})|<r$ if $P(s)$ is small and $L(\lambda, \lambda, s)$ is not very small on the edge of the disc. Thus we need a bound from below for $L(\lambda, \lambda, s)$ when $s$ is close to a zero. 
\begin{Proposition}\label{1/L}
Let $0<\lambda, \alpha\le1$. Let $\sigma_0\in\mathbb R$ and $\Re s\ge\sigma_0$. Let $L(\lambda,\alpha,s)\ne0$ and  $d$ be the distance from $s$ to the nearest zero of $L(\lambda,\alpha,s)$. Then
\begin{align*}
\frac1{|L(\lambda,\alpha, s)|}<\exp(C(|\log d|+1)\log t),
\end{align*}
where $C=C(\lambda, \alpha, \sigma_0)$ is a positive constant.
\end{Proposition}
 The proposition will help us to prove the following theorem.
\begin{Theorem}\label{localalmostsymmetry}
Let $0<\lambda\le1$ and $A>0$ be such that $AC<\pi$, where $C=C(\lambda,\alpha, -1)$ is from Proposition \ref{1/L}. Let $\rho=\beta+i\gamma$ be a nontrivial zero of $L(\lambda, \lambda, s)$. If $\gamma$ is sufficiently large, then there is a radius $r$, 
\begin{align}\label{inqexp}
\exp(-A\gamma/\log \gamma)\le r\le \exp(-A\gamma/\log \gamma)\log^2\gamma,
\end{align}
such that 
the discs
$$|s-\rho|<r\quad\text{and}\quad |s-(1-\overline{\rho})|<r$$
contain the same number of zeros.
\end{Theorem}

  In the next section we present the computer calculations related to Theorem \ref{localalmostsymmetry}. Sections \ref{proofTh1},  \ref{proofPr2}, and \ref{proofTh3} contain proofs of Theorem \ref{averagesymm}, Proposition \ref{1/L}, and Theorem \ref{localalmostsymmetry} respectively.

\section{Computations}

This section  is devoted to the more precise calculations of the first nontrivial zeros. If a nontrivial zero $\rho$ of $L(\lambda, \lambda, s)$ lies on the critical line, then by the functional equation \eqref{almostsymetry} we have that $L(1-\lambda, 1-\lambda, \rho)=0$. Similarly, if $L(\lambda, \lambda, s)$ has symmetrical zeros $\rho$ and   $1-\overline{\rho}$ then again  $L(1-\lambda, 1-\lambda, \rho)=0$. Let  $\rho_1=0.50...+9.69... i$, $\rho_2=0.50...+15.26... i$, $\rho_3=0.50...+18.65... i$, $\rho_4=0.50...+23.05... i$ be the first four zeros of $L(3/4,3/4, s)$ indicated in the Introduction. We have that

 $|L(1/4, 1/4, \rho_1)|=2.73...$,
 
 $|L(1/4, 1/4, \rho_2)|=0.13...$, 
  
 $|L(1/4, 1/4, \rho_3)|=0.48...$, 
 
 $|L(1/4, 1/4, \rho_4)|=1.15...$.  
 
 \noindent
  Thus the zeros  $\rho_1, \rho_2, \rho_3, \rho_4$ of $L(3/4,3/4, s)$ do not lie on the critical line. Using arbitrary-precision floating-point arithmetic computations, we get that 
  
  $\Re\rho_1=0.5+7.16...\cdot10^{-14}$,
  
  $\Re\rho_2=0.5-6.08...\cdot10^{-23}$,
  
  $\Re\rho_3=0.5-4.53...\cdot10^{-27}$,
  
  $\Re\rho_4=0.5-1.11...\cdot10^{-32}$. 

\noindent
The last four lines were computed in the following two ways: one by using {\it findroot} and the other by computing the contour integral which encloses only one zero $\rho$ of $L(3/4,3/4, s)$. For more details on computation methodology see the end of this section.
 
 In the upper half-plane, the first pair of almost symmetrical zeros of $L(3/4,3/4, s)$  is $-0.10...+120.59...i$ and $1.10...+120.59...i$. These zeros are not strictly symmetrical, since
 
  $|L(1/4, 1/4, 1.10...+120.59...i)|=3.94...\ne0,$ 
  
  $|L(1/4, 1/4,-0.10...+120.59...i)|=23.49...\ne0.$
  
 \noindent  
  Further, we give a table (see Table \ref{tab:number_of_zeros}) where the number of nontrivial zeros in $0<t<300$ is calculated for various cases of $L(\lambda, \lambda, s)$. For all those zeros, we have checked that $L\left(\lambda,\lambda,\rho\right) = 0$ implies $L\left(1-\lambda,1-\lambda,\rho\right) \neq 0$ if $\lambda \neq \frac{1}{2}$ . Thus, in Table \ref{tab:number_of_zeros}, all zeros, except the case $\lambda=1/2$,  are not strictly symmetrical with respect to the critical line.  

\begin{table}
\caption{Distribution of   the nontrivial zeros of $L\left(\lambda,\lambda,\sigma+it\right)$ in $0<t<300$. $N_1$ is the total number of the nontrivial zeros; $N_2$ is the number of the nontrivial zeros satisfying $|\Re \rho-1/2|>10^{-9}$, these zeros appear in almost symmetrical pairs; in the last column of the table, we have $100N_2/N_1$. } \label{tab:number_of_zeros} 
\begin{center}
\begin{tabular}{l|r|r|r}
$\lambda$ & $N_1$ & $N_2$ &  \%  \\
\hline
1/2 &              203 &                             0 &                          0.00 \\
5/9 &              193 &                            28 &                         14.51 \\
4/7 &              191 &                            24 &                         12.57 \\
3/5 &              186 &                            14 &                          7.53 \\
5/8 &              182 &                            22 &                         12.09 \\
2/3 &              176 &                            18 &                         10.23 \\
7/10 &              171 &                            28 &                         16.37 \\
5/7 &              169 &                            30 &                         17.75 \\
3/4 &              165 &                            20 &                         12.12 \\
7/9 &              161 &                            26 &                         16.15 \\
4/5 &              159 &                            22 &                         13.84 \\
5/6 &              155 &                            22 &                         14.19 \\
6/7 &              151 &                            28 &                         18.54 \\
7/8 &              150 &                            30 &                         20.00 \\
8/9 &              149 &                            22 &                         14.77 \\
9/10 &              147 &                            24 &                         16.33 \\
\end{tabular}
\end{center}
\end{table}

Computations were validated with the help of Python with {\it mpmath}\footnote{Fredrik Johansson and others. mpmath: a Python library for arbitrary-precision floating-point arithmetic (version 0.18), December 2013. http://mpmath.org/.} package.  We used the following expression of the Lerch zeta-function for rational parameters  
\begin{align*}
L\left(s, \frac{b}{d},\frac{b}{d}\right) & =\sum_{k=0}^{d-1}\sum_{m=0}^{\infty}\frac{\exp\left(2\pi i\frac{b}{d}\left(dm+k\right)\right)}{\left(dm+k+\frac{b}{d}\right)^{s}}\\
& =d^{-s}\sum_{k=0}^{d-1}\exp\left(2\pi i\frac{b}{d}k\right)\zeta\left(s,\frac{kd+b}{d^{2}}\right),
\end{align*}
where $\zeta(s,\alpha)$, $0<\alpha\le1$, is the Hurwitz zeta-function.  The function $\zeta(s,\alpha)$ is
implemented by  the command {\it zeta}.  Zero locations were calculated using {\it findroot} with Muller\textquoteright s method. 

In this paper, all computer computations  should be regarded as heuristic because their accuracy was not controlled explicitly.

\section{Proof of Theorem \ref{averagesymm}}\label{proofTh1}

In \cite{gs2002}, it was proved that, for  $0<\lambda,\alpha\leq 1$,
\begin{align}\label{modulus}
\sum_{|\gamma|\leq T}\left(\beta-{1\over 2}\right)
={T\over 2\pi}\log{\alpha\over \sqrt{\lambda(1-\{\lambda\})}}+O(\log T).
\end{align}
Theorem \ref{averagesymm} can be derived from the proof of formula (\ref{modulus}).  Namely, from the proof of Theorem 1 in \cite{gs2002} we derive the following lemma.
\begin{Lemma}
Let $b\geq 3$ be a constant. For $0<\lambda,\alpha\leq 1$,
$$ 
\sum_{0<\gamma\leq T}(b+\beta)=\left(b+{1\over 2}\right)
{T\over 2\pi}\log{T\over 2\pi e\alpha\lambda}+{T\over
4\pi}\log{\alpha\over\lambda}
+O(\log T).$$
\end{Lemma}
Then the equality
\begin{align*}
\sum_{0<\gamma\leq T}\left(\beta-{1\over 2}\right)=\sum_{0<\gamma\leq T}(b+\beta)-(b+\frac12)\sum_{0<\gamma\leq T}1
\end{align*}
together with the zero counting formula (\ref{zeronumber})   gives Theorem~\ref{averagesymm}.

\section{Proof of Proposition \ref{1/L}}\label{proofPr2}

We start from the following lemma.
\begin{Lemma}\label{fbyzeros}
If $f(s)$ is regular, and
\begin{eqnarray*}
\left\vert{f(s)\over f(s_0)}\right\vert<e^M
\end{eqnarray*}
in $\{s\,:\,\vert s-s_0\vert\leq r\}$ with $M>1$, then
\begin{eqnarray*}
\left\vert{f(s_0)\over f(s)}\prod_\rho{s-\rho\over
s_0-\rho}\right\vert<e^{CM}
\end{eqnarray*}
for $\vert s-s_0\vert\leq {3\over 8}r$, where $C$ is some constant
and $\rho$ runs through the zeros of $f(s)$ such that
$\vert\rho-s_0\vert\leq {1\over 2}r$.
\end{Lemma}
\proof 
The lemma  follows immediately from the proof  of Lemma $\alpha$ in Titchmarsh \cite[\S 3.9]{Titchmarsh1986}.
\endproof

In order to apply Lemma \ref{fbyzeros}, we need information about the growth of the Lerch zeta-function. For each $\sigma$, we define a number $\mu(\sigma)=\mu(\lambda, \alpha, \sigma)$ as the lower bound of numbers $\xi$ such that $L(\lambda,\alpha,\sigma+iT)\ll T^\xi$.
\begin{Lemma}\label{lerchgrow}
Let $0<\lambda,\alpha\le1$ and $\sigma_0<0$. Then
\begin{align*}
\mu(\sigma)\le
\begin{cases}
\frac12-\sigma &\text{if } \sigma_0\le\sigma\le0, 
\\
\frac12+\left(\frac{64}{205}-1\right)\sigma &\text{if } 0\le\sigma\le\frac12, 
\\
\frac{64}{205}(1-\sigma) &\text{if } \frac12\le\sigma\le1, 
\\
0 &\text{if } \sigma\ge1.
\end{cases}
\end{align*}
\end{Lemma}
\proof
 In \cite{Garunkstis2005}, it is proved that 
$$L(\lambda,\alpha,\frac12+it)\ll t^{\frac{32}{205}+\varepsilon}\qquad(t\to\infty),$$
and from the approximation of the Lerch zeta-function by a finite sum  (\cite[Theorem 1.2 in Chapter 3]{gl}) we see that $L(\lambda,\alpha,\sigma+it)\ll t^{\varepsilon}$ for $\sigma\ge1$.
Now the lemma follows by the Phragm\'en-Lindel\"of theorem (see Titchmarsh \cite[\S 5.65]{tit}) and by  the functional equation (\ref{Lerchfunc}) in view of Stirling's formula (see Titchmarsh \cite[\S 4.42]{tit})
\begin{align}\label{gamma}
|\Gamma(1-s)|=\sqrt{2\pi}|t|^{\frac12-\sigma}
e^{-\frac{\pi |t|}2}(1+O(|t|^{-1}))\qquad (|t|\to\infty), 
\end{align}
uniformly  for $\sigma_0<\sigma\le1/2$.
\endproof

\proof[Proof of Proposition \ref{1/L}.] To prove the proposition  we choose $f(s)=L(\lambda,\alpha, s)$, $s_0=3+it$, and a sufficiently large but fixed radius $r$ in Lemma \ref{fbyzeros}. In view of Lemma \ref{lerchgrow} we take  $M=b\log T$, where $b=b(r)$. The function $1/L(\lambda,\alpha, s_0)$ is bounded. By the formula  (\ref{zeronumber}) for a number of nontrivial zeros, we have that the number of zeros in the disc $|s-s_0|<\frac12r$ is $<c\log \Im s_0$. This proves Proposition \ref{1/L}.

\endproof

\section{Proof of Theorem \ref{localalmostsymmetry}}\label{proofTh3}

\proof[Proof of Theorem \ref{localalmostsymmetry}]

If  $\lambda=1/2,1$, then in view of equalities \eqref{RH} it is well known that the non-real complex number $\rho$ is a zero of $L(\lambda,\lambda, s)$ if and only if $1-\overline{\rho}$ is also a zero of $L(\lambda,\lambda, s)$.

Next we assume that $0<\lambda<1$ and $\lambda\ne1/2$. By the formula  (\ref{zeronumber}), we see that the number of zeros in the disc $|s-\rho|<\exp(-A\gamma/ \log^2\gamma)$ is $<c\log \Im\rho$. Let $r_k=k\exp\left(-A\gamma/\log \gamma\right)$, $k=1,\dots,[c\log \gamma]+1$. By Dirichlet's box principle, there is $1\le\ell\le[c\log \gamma]$ such that $L(\lambda,\lambda, s)$ has no zeros for the ring
\begin{align*}
r_\ell<|s-\rho|\le r_{\ell+1}.
\end{align*}
Let
\begin{align*}
r=(r_\ell+r_{\ell+1})/2=(\ell+1/2)\exp\left(-A\gamma/\log \gamma\right).
\end{align*}

 Suppose $G(s)$ and $P(s)$ are defined by the functional equation (\ref{almostsymetry}). Let  $C_R=\{s : |s-\rho|=r\}$ and $C_L=\{s : |s-(1-\overline{\rho})|=r\}$.   The steps of the proof are the following.  If $L(\lambda, \lambda, s)$ has  $N$ zeros inside of $C_R$, then  by Rouch\'e's theorem we expect that $G(s)L(\lambda, \lambda, s)+P(s)$ has $N$ zeros inside $C_R$. Then by the functional equation (\ref{almostsymetry}), the function   $L(1-\lambda,\lambda,1-s)$ has $N$ zeros inside $C_R$, then by conjugation $L(\lambda, \lambda, s)$  has $N$ zeros inside $C_L$.  Next we need to justify the step involving  Rouch\'e's theorem.

Note that $G(s)$ has no zeros. By Rouch\'e's theorem, the functions $G(s)L(\lambda,\lambda,s)$ and $G(s)L(\lambda,\lambda,s)+P(s)$ have the same number of zeros inside of the circle $C_r$ if on this circle the inequality 
\begin{align}\label{P(s)<}
|P(s)|<|G(s)L(\lambda,\lambda,s)|
\end{align} 
is valid. 

 In view of the growth of the Lerch zeta-function (see Lemma \ref{lerchgrow})  we get that, for sufficiently large $t$ and $-1.4\le\sigma\le2$, 
$$|P(s)|< |\Gamma(s)|t^2e^{-\pi t/2}\quad \text{and}\quad |G(s)|\ge (2\pi)^{-2} |\Gamma(s)|e^{\pi t/2}.$$ 
 Proposition \ref{1/L} gives, for $s\in C_R$, 
 \begin{align*}
 L(\lambda, \lambda, s)\gg\exp\left(-\left(AC + o(1)\right) \gamma\right),
 \end{align*}
where $AC<\pi$.
Thus the inequality (\ref{P(s)<}) is valid.   
By this Theorem \ref{localalmostsymmetry} is proved.

\endproof

Note, that from this proof we have that the quantity $\log^2\gamma$ in the inequality \eqref{localalmostsymmetry} of Theorem \ref{localalmostsymmetry}  can be replaced (at the expense of more complicated notations)  by the smaller quantity $c\log\gamma+1$, where $c$ is from the proof of Theorem \ref{localalmostsymmetry}.


\begin{thebibliography}{99}

\bibitem{Apostol1951}
{\sc T. M. Apostol}, {\em  On the Lerch zeta function},
 Pacific J. Math. {\bf 1} (1951), 161--167.  {\em Addendum}, Pacific J. Math. {\bf 2}(1952), 10. 

\bibitem{Balanzario2007}
{\sc E. P. Balanzario, J. Sanchez-Ortiz}, {\em Zeros of the Davenport-Heilbronn counterexample},
Math. Comput. {\bf 76} (2007), 2045--2049.

\bibitem{Berndt1972}
{\sc B. C. Berndt}, {\em Two new proofs of Lerch's functional equation},
Proc. Amer. Math. Soc. {\bf 32} (1972), 403--408.

\bibitem{bh}
{\sc E. Bombieri, D. A. Hejhal}, {\em On the distribution of zeros of linear combinations of Euler products},
Duke Math. J. {\bf 80} (1995), 821--862.

 \bibitem{Garunkstis2005}
{\sc R. Garunk\v stis}, {\em Growth of the Lerch zeta-function},
  Lith. Math. J. {\bf 45} (2005), 34--43. 


\bibitem{gar}{\sc R. Garunk\v{s}tis, A. Laurin\v{c}ikas}, {\it On
zeros of the Lerch zeta-function}, Number theory and its
applications, S. Kanemitsu and K. Gy\"ory (editors), Kluwer
Academic Publishers 1999, 129--143.


\bibitem{gls2003}{\sc R. Garunk\v{s}tis, R. A. Laurin\v{c}ikas, J. Steuding}, {\it On the mean square of Lerch zeta-functions},  Arch. Math. {\bf 80} (2003), 47--60.

\bibitem{gs2002}{\sc R. Garunk\v{s}tis, J. Steuding}, {\it On the zero distributions of Lerch zeta-functions}, Analysis {\bf 22} (2002), 1--12.

\bibitem{Garunkstis2002}
{\sc R. Garunk\v stis, J. Steuding}, {\em Do Lerch zeta-functions satisfy the Lindelof hypothesis?}, in: Analytic and Probabilistic Methods in Number Theory, Proceedings of the Third Intern. Conf. in Honour of J. Kubilius, Palanga, Lithuania, 24-28 September 2001, (eds. A. Dubickas, A. Laurin\v cikas and E. Manstavi\v cius), TEV, Vilnius, (2002), 61--74.


\bibitem{gs}{\sc R. Garunk\v{s}tis, R. \v Sim\.enas}, {\it On the Speiser equivalent for the Riemann hypothesis},  Eur. J. Math. {\bf 1} (2015), 337--350.


\bibitem{ll1}{\sc J. C. Lagarias, W.-C. W. Li}, {\it The Lerch zeta function I. Zeta integrals},  Forum Math. {\bf 24} (2012), 1--48.

\bibitem{ll2}{\sc J. C. Lagarias, W.-C. W. Li}, {\it The Lerch zeta function II. Analytic continuation},  Forum Math. {\bf 24} (2012), 49--84.

\bibitem{lau}{\sc A. Laurin\v cikas}, {\it Universality of the Lerch zeta function},  Lithuanian Math. J.  {\bf 37} (1997), 275--280.

\bibitem{gl}{\sc A. Laurin\v cikas, R. Garunk\v stis}, {\it The Lerch zeta-function},
 Kluwer Academic Publishers, 2002, 197 pp.
 
 \bibitem{lnp}
{\sc Y. Lee, T. Nakamura, \L. Pa\'nkowski}, {\em Joint universality for Lerch zeta-functions},
J. Math. Soc. Japan {\bf 69} (2017), 153--161.

\bibitem{Lerch1887}
{\sc M. Lerch}, {\em Note sur la fonction ${\mathcal K}(z,x,s)=
\sum_{k=0}^\infty e^{2k\pi ix}(z+k)^{-s}$},
Acta Math. {\bf 11}
(1887), 19--24.

\bibitem{Mikolas1971}
{\sc M. Mikol\'as}, {\em New proof and extension of the functional equality of Lerch's zeta-function},
Ann. Univ. Sci. Budapest. Sect. Math. {\bf 14} (1971), 111--116.

\bibitem{Oberhettinger1956}
{\sc F. Oberhettinger}, {\em Note on the Lerch zeta-function},
 Pacific J. Math. {\bf 6} (1956), 117--120. 

\bibitem{Spira}
{\sc R. Spira}, {\em Zeros of Hurwitz zeta-functions},
 Math Comput. {\bf 136} (1976), 863--866.

\bibitem{tit}{\sc E. C. Titchmarsh}, {\it The theory of functions}, 2nd ed.,
Oxford University Press, 1939.

\bibitem{Titchmarsh1986}
{\sc E. C. Titchmarsh}, {\em The theory of the Riemann zeta-function. 2nd ed., rev. by D. R. Heath-Brown}, Oxford Science Publications. Oxford: Clarendon Press, 1986.

\bibitem{vaug}
{\sc R. C. Vaughan}, {\em Zeros of Dirichlet series},
Indagationes Mathematicae {\bf 26} (2015), 897--909.

\end{thebibliography}
\end{document}